\documentclass{article}

\usepackage[utf8]{inputenc}
\usepackage[cmintegrals]{newtxmath}
\usepackage[T1]{fontenc}
\usepackage[a4paper, total={6in, 8in}]{geometry}
\usepackage{graphicx}
\usepackage{mathrsfs}
\usepackage{amsmath,stackrel}
\usepackage{graphicx}
\usepackage{multirow}
\usepackage{mathtools}
\usepackage{ marvosym }
\usepackage{tikz-cd}
\usepackage{pst-node}

\newcommand{\so}{\mathfrak{so}(3)}

\newcommand{\N}{\mathbb{N}}

\newcommand{\R}{\mathbb{R}}
\newcommand{\grad}{\text{grad}}

\newcommand{\J}{\text{J}}

\newcommand{\SO}{\text{SO}}

\usepackage{bm}
\usepackage[colorlinks=true]{hyperref}
\hypersetup{urlcolor=blue, citecolor=red}
\newtheorem{theorem}{Theorem}[section]
\newtheorem{corollary}{Corollary}
\newtheorem{lemma}[theorem]{Lemma}
\newtheorem{proposition}{Proposition}

\newtheorem{definition}[theorem]{Definition}
\newtheorem{remark}{Remark}

\title{Variational Collision Avoidance on Riemannian manifolds}
\author{Jacob R. Goodman and Leonardo J. Colombo}
\date{}

\begin{document}

\maketitle

\begin{abstract}
This paper studies variational collision avoidance problems for multi-agents systems on complete Riemannian manifolds. That is, we minimize an energy functional, among a set of admissible curves,which depends on an artificial potential function used to avoid collision between the agents. We show the global existence of minimizers to the variational problem and we provide conditions under which it is possible to ensure that agents will avoid collision within some desired tolerance. We also study the problem where trajectories are constrained to have uniform bounds on the derivatives, and derive alternate safety conditions for collision avoidance in terms of these bounds - even in the case where the artificial potential is not sufficiently regular to ensure existence of global minimizers. 
\end{abstract}

\section{Introduction}
Energy-optimal path planning on nonlinear spaces such as Riemannian manifolds has been an active field of interest in the last decades due to its numerous applications in manufacturing, aerospace technologies, and robotics \cite{blochgupta, bonard, hussein, ring, tar, trouve, zefran}. It is often the case that the desired paths must connect some set of knot points—interpolating positions with given velocities and potentially higher order derivatives \cite{CLACC, CroSil:95, mel}. For such problems, the use of variationally defined curves has a rich history due to the regularity and optimal nature of the solutions. In particular, the so-called Riemannian splines \cite{noakes} are a particularly ubiquitous choice in interpolant, which themselves are composed of Riemannian polynomials—satisfying boundary conditions in positions, velocities, and potentially higher-order derivatives—that are glued together. In Euclidean spaces, Riemannian splines are just cubic splines, that is, the minimizers of the total squared acceleration \cite{boor}.

Riemannian polynomials are smooth and optimal in the sense that they minimize the average square magnitude of some higher-order derivative (a quantity which is often related to energy consumption in applications). Moreover, Riemannian polynomials carry a rich geometry with them, which has been studied extensively in the literature (see \cite{Giambo, marg, noakes, elastica} for a detailed account of Riemannian cubics and \cite{RiemannianPoly}, \cite{popei} for some results with higher-order Riemannian polynomials).

It is often the case that—in addition to interpolating points—there are obstacles or regions in space which need to be avoided. In this case, a typical strategy is to augment the action functional with an artificial potential term that grows large near the obstacles and small away from them (in that sense, the minimizers are expected to avoid the obstacles) \cite{kod}, \cite{rybus}, \cite{chang}. This was done for instance in \cite{BlCaCoCDC} and \cite{BlCaCoIJC}, where necessary conditions for extrema in obstacle avoidance problems on Riemannian manifolds were derived, in addition to applications to interpolation problems on manifolds and to energy-minimum problems on Lie groups and symmetric spaces endowed with a bi-invariant metric. Energy-minimum obstacle avoidance problems on Riemannian manifolds for nonholonomic systems were studied in \cite{colombononholonomic}. More recent approaches based on the use of artificial potentials include hybrid feedback controllers \cite{piveda, ricardo1, ricardo2, dimoshybrid} and source seeking control based-approaches \cite{dur}, \cite{mohr}, among others. 

In the case of path planning for multi-agent systems, another practical consideration is often necessary. Namely, that the agents do not collide along their trajectories in addition to the individual tasks (such as interpolation, obstacle avoidance, and minimizing some cost functional). This can be handled with a similar strategy to that of obstacle avoidance. That is, by augmenting the action to be minimized with an artificial potential that grows large when two agents, which can sense each other, are sufficiently close —as was done in \cite{sh}, \cite{CoGo20} for  Riemannian manifolds and \cite{dimos2} for Euclidean spaces (see also \cite{piveda} and \cite{gram} for other recent approaches based on hybrid feedback stabilization on Euclidean spaces). An important point to consider when applying such a methodology to applications—which thus far has been lacking in the literature—is that of safety guarantees. That is, can an artificial potential be designed to ensure that agents will avoid collision within some desired tolerance, and if so, then how. Answering these questions is the main focus of this paper. 

In particular, our aim is to rigorously investigate the role of the artificial potential in the variational collision avoidance task on complete and connected Riemannian manifolds, and in doing so obtain some conditions under which safety is guaranteed. The three main contributions of this paper are as follows. (i) We prove the existence of global minimizers to the variational problem in the case that the potential is $C^1$ and non-negative, which is a necessary prerequisite in providing safety guarantees (indeed, proving that minimizing trajectories avoid collision is useful only if such minimizing trajectories exist). This is accomplished by showing that the functional in which we are minimizing satisfies the Palais-Smale condition on its domain, as was done in \cite{Giambo}, \cite{RiemannianPoly} for Riemannian polynomials in the single-agent setting. (ii) We derive general conditions for the artificial potential—in terms of some reference trajectory which avoids collision—under which the corresponding minimizers avoid collision within some tolerance. We then remove the dependence on the reference trajectory for a particular family of potentials, and show that collision can be avoided within any desired tolerance for some potential in the family (constrained by the boundary conditions and the geometry of the manifold). (iii) We study the problem where trajectories are constrained to have uniform bounds on the derivatives, and derive alternate safety conditions for collision avoidance in terms of these bounds (even in the case where the artificial potential is not sufficiently regular to ensure existence of global minimizers). 

The remainder of the paper is structured as follows. Section \ref{Sec: background} provides some background in Riemannian geometry and Sobolev spaces of curves on Riemannian manifolds — which serve as the natural domain of our action functional. In Section \ref{Sec: Necessary conditions}, we define the variational problem that we wish to solve, and provide necessary conditions for optimality in Proposition \ref{th1}. Section \ref{Sec: existence} regards the existence of global minimizers, and it is proven in Theorem \ref{PS} that our functional satisfies the Palais-Smale condition on its domain. Section \ref{Sec: Collision Avoidance} contains the main results regarding safety guarantees for collision avoidance. In particular, we define the collision avoidance task and derive conditions under which collision avoidance is guaranteed in Proposition \ref{Collision_Avoidance_Prop} and Proposition \ref{Prop: Coll_avoid_bounded}. We end up the paper by providing some simulation results to show how the main results of this paper can be applied in concrete situations.

\section{Sobolev Spaces of Curves on a Riemannian Manifold}\label{Sec: background}

Let $Q$ be an $n$-dimensional  \textit{Riemannian
manifold} endowed with a symmetric covariant 2-tensor field $g$ called the \textit{Riemannian metric}. That is, to each point $q\in Q$ we assign an inner product $g_q:T_qQ\times T_qQ\to\mathbb{R}$, where $T_qQ$ is the \textit{tangent space} of $Q$ at $q$. The length of a tangent vector is determined by its norm,
$||v_q||=g(v_q,v_q)^{1/2}$ with $v_q\in T_qQ$. A \textit{Riemannian connection} $\nabla$ on $Q$ is a map that assigns to any two smooth vector fields $X$ and $Y$ on $Q$ a new vector field, $\nabla_{X}Y$. For the properties of $\nabla$, we refer the reader to \cite{Boothby,Milnor}.  The operator
$\nabla_{X}$, which assigns to every vector field $Y$ the vector
field $\nabla_{X}Y$, is called the \textit{covariant derivative of
$Y$ with respect to $X$}. Given $V:Q\to\mathbb{R}$, we define the gradient vector field,  $\hbox{grad}_{g}V$ as the vector field on $Q$ characterized by $g(\hbox{grad}_{g}V, X) =
X(V), \mbox{ for  every vector field } X$ on $Q$.


Consider a vector field $W$  along a curve $q$ on $Q$. The $k$th-order covariant derivative  of $W$ along $q$ is denoted by $\displaystyle{\frac{D^{k}W}{dt^{k}}}$, $k\geq 1$. We also denote by $\displaystyle{\frac{D^{k+1}q}{dt^{k+1}}}$ the $k$th-order covariant derivative of the velocity vector field of $q$  along $q$, $k\geq 1$.

A vector field $X$ along a piecewise smooth curve $q$ in $Q$ is said to be \textit{parallel along $q$} if $\displaystyle{\frac{DX}{dt}\equiv 0}$.

Given vector fields $X$, $Y$
and $Z$ on $Q$, the vector field $R(X,Y)Z$ given by \begin{equation}\label{eq:CurvatureEndomorphismDefinition}
R(X,Y)Z=\nabla_{X}\nabla_{Y}Z-\nabla_{Y}\nabla_{X}Z-\nabla_{[X,Y]}Z
\end{equation}  is called the \textit{curvature endomorphism} on $Q$. $R$ is trilinear in $X$, $Y$ and $Z$. We further define the \textit{curvature tensor} on $Q$ by \begin{equation}\label{eq:CurvatureTensorDefinition}
Rm(X, Y, Z, W) = g\left(R(X, Y)Z, W\right).
\end{equation}

Let $\Omega $ be the set of all $\mathcal{C}^1$ piecewise smooth curves $q:[0,T] \to Q$ such that $q(0)$, $q(T)$, $\frac{dq}{dt}(0)$, $\frac{dq}{dt}(T)$ are fixed. The set $\Omega$ is called the \textit{admissible set}. A $\mathcal{C}^1$-piecewise smooth one-parameter \textit{admissible variation} of a curve $q\in\Omega$ is a family of curves $\alpha:(-\epsilon,\epsilon)\times [0,T]\to Q$; $(r,t)\to \alpha(r,t)=\alpha_r(t)$ such that $\alpha_0=q$ and $\alpha_r\in\Omega$ for each $r\in (-\epsilon,\epsilon)$. 

If we assume that $Q$ is \textit{complete}, then any two points $x$ and $y$ in $Q$ can be connected by a geodesic, and the Riemannian distance between two points in $Q$, $d:Q\times Q\to\mathbb{R}$ can be defined by
$\displaystyle{d^2(q,y)=\int_{0}^{1}\Big{\|}\frac{d \gamma_{q,y}}{d s}(s)\Big{\|}^2\, ds}$, where $ \gamma_{q,y}$ is the geodesic connecting the points $q$ and $y$ on $Q$. The idea of a geodesic is useful because it provides a map $T_qQ\to Q$ in the following way
\begin{equation*}
    v \mapsto \gamma(1),\quad \gamma(0)=q,\quad \dot{\gamma}(0)=v
\end{equation*}
where $\gamma$ is a geodesic. This map is called the Riemannian exponential map and is denoted by $\mathrm{exp}_q:T_qQ\to Q$. In particular, $\mathrm{exp}_q$ is a diffeomorphism from some star-shaped neighborhood of $0 \in T_q Q$ to a convex open neighborhood $\mathcal{B}$ of $q \in Q$, and if $y \in \mathcal{B}$, we can write the Riemannian distance by means of the Riemannian exponential as $d(q,y)=\|\mbox{exp}_q^{-1}y\|.$

\vspace{.2cm}

The \textit{Lebesgue space} $L^p([0,1];\R^n)$, $p\in(1,+\infty)$ is the space of $\R^n$-valued functions on $[0,1]$ such that each of their components is $p$-integrable, that is, whose integral of the absolute value raised to the power of $p$ is finite. 

A sequence $(f_n)$ of functions in $L^p([0,1];\R^n)$ is said to be \textit{weakly convergent} to $f$ if for every $g\in L^r([0,1];\R^n)$, with  $\frac{1}{p}+\frac{1}{r}=1$, and every component $i$,  $\displaystyle{\lim_{n\to\infty}\int_{[0,1]}f_n^i g^i=\int_{[0,1]}f^i g^i}$.

A function $g\colon[0,1]\to \R^n$ is said to be the \textit{weak derivative} of $f\colon[0,1]\to \R^n$ if for every component $i$ of $f$ and $g$, and for every compactly supported $\mathcal{C}^\infty$ real-valued function $\varphi$ on $[0,1]$, $\displaystyle{\int_{[0,1]}f^i\varphi'=-\int_{[0,1]}g^i\varphi}$.

The \textit{Sobolev space} $W^{k,p}([0,1];\R^n)$ is the space of functions $u\in L^p([0,1];\R^n)$ such that for every $\alpha\leq k$, the $\alpha^{th}$ weak derivative $\frac{d^\alpha u}{dt^{\alpha}}$ of $u$ exists and $\frac{d^\alpha u}{dt^{\alpha}}\in L^p([0,1];\R^n)$. In particular, $H^k([0,1];\R^n)$ denotes the Sobolev space $W^{k,2}([0,1];\R^n)$. A sequence $(f_n)\subset W^{k,p}([0,1];\R^n)$ is said to be \textit{weakly convergent} to $f$ in $W^{k,p}([0,1];\R^n)$ if for every $\alpha\leq k$, $\displaystyle{\frac{d^\alpha f_n}{dt^{\alpha}}\rightharpoonup \frac{d^\alpha f}{dt^{\alpha}}}$ weakly in $L^p([0,1];\R^n)$.


Let $Q$ be an $m$-dimensional Riemannian manifold. By $H^2([0,1];Q)$ we will denote the set of all curves $q\colon[0,1]\to Q$ such that for every chart $(\mathcal{U},\varphi)$ of $Q$ and every closed subinterval $I\subset[0,1]$ such that $q(I)\subset\mathcal{U}$, the restriction of the composition $\varphi\circ q|_I$ is in $H^2([0,1];\R^m)$. Note that $H^2([0,1]; Q)$ is an infinite-dimensional Hilbert Manifold modeled on $H^2([0,1]; \R^m)$, and given $\xi = (p,v)\in T_p Q$ and $\eta = (q, w) \in T_q Q$, the space $$\Omega_{\xi, \eta} := \{\gamma \in H^2([0,1]; Q) : \gamma(0) = p, \ \gamma(T), \ \dot{\gamma}(0) = v, \ \dot{\gamma}(T) = w \}$$ is a closed submanifold of $H^2([0,1]; Q)$ (see \cite{Palais}, \cite{Palais2}, \cite{Piccione} for instance). The tangent space $T_x \Omega_{\xi, \eta}$ consists of vector fields along $x$ of class $H^2$ which vanish at the endpoints together with their first covariant derivatives. We consider the Hilbert structure on $T_x \Omega_{\xi, \eta}$ induced by the inner product: $$\left< V, W \right> := \int_0^T g\left(\frac{D^2}{\partial t^2}V, \frac{D^2}{\partial t^2}W\right)dt.$$
Similarly, we may consider the product space $\Omega_{\xi_1, \eta_1} \times \cdots \times \Omega_{\xi_s, \eta_s}$ as a closed submanifold of the Hilbert manifold consisting of $s$ copies of $H^2([0,1]; Q)$, endowed with the inner product on $T_{(x_1,...,x_s)} (\Omega_{\xi_1, \eta_1} \times \cdots \times \Omega_{\xi_s, \eta_s}) \cong T_{x_1} \Omega_{\xi_1, \eta_1} \times \cdots \times T_{x_s} \Omega_{\xi_s, \eta_s}$:

$$\left< (V_1,..., V_s), (W_1,..., W_s)\right> := \sum_{i=1}^s \left<V_i, W_i\right>.$$


\section{Variational Collision Avoidance for Multi-agent Systems}\label{Sec: Necessary conditions}
Consider a set $\mathcal{V}$ consisting of $s\geq 2$ agents on the Riemannian manifold $Q$. The configuration of each agent at any given time is determined by the element $q_i\in Q$, $i=1,\ldots,s$. The neighboring relationships are described by an undirected time-invariant graph $\mathbb{G} = (\mathcal{V}, \mathcal{E})$ with edge set $\mathcal{E}\subseteq\mathcal{V}\times\mathcal{V}$. The set of neighbors $\mathcal{N}_i$ for the agent $i\in\mathcal{V}$ is given by $\mathcal{N}_i=\{j\in\mathcal{V}:(i,j)\in\mathcal{E}\}$. An agent $i\in\mathcal{V}$ can measure its Riemannian distance from other agents in the subset $\mathcal{N}_i \subseteq \mathcal{V}$.


For $i=1,...,s$, consider the set $\Omega_i$ of all $\mathcal{C}^1$-piecewise smooth curves on $Q$, $q_i:[0,T]\to Q$ satisfying the boundary conditions
\begin{equation}\label{ic}
q_i(0)=q_i^0, \;\;\; \frac{dq_i}{dt}(0)={v_i^0}, \quad
q_i(T)=q_i^T, \;\;\;\; \frac{dq_i}{dt}(T)=v_i^T,
\end{equation} 
and the functional $J$ on $\Omega=\Omega_1\times\cdot\cdot\cdot\times\Omega_s$,
\begin{equation}\label{J}
J(q_1,q_2,\ldots,q_s)=\frac{1}{2} \sum\limits_{i=1}^s \int\limits_0^T \Big{(}\Big{|}\Big{|}\frac{D^2q_i}{dt^2}(t)\Big{|}\Big{|}^2 + \frac{1}{2}\sum\limits_{j\in \mathcal{N}_i}V_{ij}(q_i(t),q_j(t))\Big{)}dt.
\end{equation}

\textit{Problem:} Find a collection of curves $(q_1,\ldots,q_s)\in \Omega$ minimizing the functional $J$ where $V_{ij}:Q\times Q\to\mathbb{R}$ is a $C^1$ non-negative artificial potential function satisfying the symmetry relations $V_{ij} = V_{ji}$ and $V_{ij}(p, q) = V_{ij}(q, p)$ for all $(i,j)\in\mathcal{E}$ and $p,q \in Q\times Q$.

\vspace{.2cm}

In order to minimize the functional $J$ among the set $\Omega$, we want to find curves $q\in\Omega$ such that $J(q)\leq J(\tilde{q})$ for all
admissible curves $\tilde{q}$ in a $C^1$-neighborhood of $q$. The next result from \cite{CoGo20} (see also \cite{mishal} for centralized communication between agents) characterizes necessary conditions for optimality in the variational collision avoidance problem.

\begin{proposition}\label{th1}\cite{CoGo20} For each $i\in\mathcal{V}$, $q_i$ must be a $\mathcal{C}^\infty$-curve on $[0,T]$ satisfying
\begin{equation}\label{eqq1}
    \frac{D^4q_i}{dt^4}+R\Big{(}\frac{D^2q_i}{dt^2},\frac{dq_i}{dt}\Big{)}\frac{dq_i}{dt}=-\sum\limits_{j\in \mathcal{N}_i} \hbox{\grad}_1 \, V_{ij}(q_i(t),q_j(t)).
\end{equation}
\end{proposition}

\begin{remark}
Note that the above formalism can be easily adapted to energy-minimum problems with different kinds of collective behavior performances other than collision avoidance including obstacle avoidance \cite{BlCaCoCDC}, consensus/synchronization on complete Riemannian manifolds \cite{LS}, \cite{MM}, interpolation among cells in a cell decomposition problem with path planning on complete manifolds \cite{TAN}, and the synchronization of quantum karumoto models on the complete manifold $SU(2)$, by employing the consensus performance given in \cite{GR}. 

Indeed, it is the shape of the potential—which we will study extensively in Section \ref{Sec: Collision Avoidance} specifically for the collision avoidance task—and the topology of the graph that decides the collective behavior. To that end, the existence of global minimizers for the functional $J$ given in equation \eqref{J} that we study in the next section will not be specific to the collision avoidance task. 
\end{remark}


\section{Existence of global minimizers}\label{Sec: existence}

Next, we show the existence of global minimizers for $J$ in $\Omega.$ In particular, we will do so by showing that $J$ satisfies the Palais-Smale condition on $\Omega$. 

We begin by defining what it means for a functional to satisfy the Palais-Smale condition, and introduce a result from $\cite{Palais}$ that motivates our consideration of the condition.

\begin{definition}
A sequence $\{q^n\} \subset \Omega$ is called a \textit{Palais-Smale sequence} for $J$ if

\vspace{.2cm}

\begin{enumerate}
    \item $\displaystyle{\sup_{n \in \N}J(q^n) < +\infty},$\\
    \item $||\mathbf{d}J(q^n)||_{T_{q^n}\Omega^\ast} \to 0$ as $n \to \infty.$
\end{enumerate}

\vspace{.3cm}

Where $\mathbf{d}J$ denotes the differential of $J$. We say that $J$ satisfies the Palais-Smale condition if every Palais-Smale sequence admits a convergent subsequence in $\Omega$.
\end{definition}

In \cite{Palais} (see Theorem $9.1.9$), it has been shown that if $f: M \to \R$ is a smooth function which is bounded below and satisfies the Palais-Smale condition on $M$—a Hilbert manifold which is complete as a Riemannian manifold—then $f$ attains its infimum in $M$. That is, there exists a critical point $x \in M$ of $f$ such that $f(x) = \displaystyle{\inf_{m \in M}f(m)}$. Given that $\Omega$ is a Hilbert Manifold, complete as a Riemannian Manifold, and $J$ is smooth and bounded below  by $0$ on $\Omega$, the existence of global minimizers follows if $J$ satisfies the Palais-Smale condition. Before proving this result, we will introduce a lemma that will simplify the analysis considerably.

\begin{lemma}\label{Convergence}
Let $Q$ be an m-dimensional complete Riemannian manifold, and suppose that $\{q^n\} = \{q_1^n,...,q_s^n\}\subset \Omega$ is a sequence such that $\displaystyle{\sup_{n \in \N} \J(q^n) < +\infty}$. Then there exists a subsequence of $\{q^n\}$ such that each $\{q_i^n\}$ converges weakly to some $q_i \in \Omega$ with respect to the $H^2$ norm.
\end{lemma}

\textit{Proof:} Suppose that $\{q^n\}$ is such a sequence. Setting $G_0 := g(v^1_0, v^1_0)$ and using the Fundamental Theorem of Calculus and the Cauchy-Shwarz inequality, we have

\begin{align*}
g\left(\dot{q}_1^n(t), \dot{q}_1^n(t)\right) &= g(\dot{q}_1^n(0), \dot{q}_1^n(0)) + \int_0^t \frac{d}{du} g\left(\dot{q}_1^n(u), \dot{q}_1^n(u)\right)du \\
& \le G_0 + 2\int_0^T \left|g\left(\frac{D^2}{\partial u^2}q_1^n(u), \dot{q}_1^n(u)\right)\right|du \\
& \le G_0 + 2 \left[ \int_0^T g\left(\frac{D^2}{\partial u^2}q_1^n(u),\frac{D^2}{\partial u^2}q_1^n(u)\right)dt \right]^{1/2} \left[ \int_0^T g\left(\dot{q}_1^n(u),\dot{q}_1^n(u)\right)du \right]^{1/2} \\
&\le G_0 + 2T \left[ \sup_{n \in \N} J(q^n) \right]^{1/2} \left[ \sup_{t \in [0,T]} g\left(\dot{q}_1^n(t),\dot{q}_1^n(t)\right) \right]^{1/2},
\end{align*}where we have used the fact that $V_{ij}$ is non-negative in the last inequality. Let $c := \displaystyle{\sup_{n \in \N}} J(q^n)$ and $G^n := \displaystyle{\sup_{t \in [0,T]}} g\left(\dot{q}_1^n(t),\dot{q}_1^n(t) \right)$. Taking the supremum of the inequality over $t \in [0, T]$, we have $$G^n \le G_0 + 2T\sqrt{cG^n} \implies G^n \le \left(T\sqrt{c} + \sqrt{T^2 c + G_0}\right)^2 := r^2.$$
Now observe that the sequence of lengths of the curves similarly satisfies $$L(q_1^n) = \int_0^T \sqrt{g\left(\dot{q}_1^n(t), \dot{q}_1^n(t)\right)}dt \le T \sqrt{G^n} \le Tr.$$

Hence, the image of $q_1^n$ is contained in the closed geodesic ball $\Bar{B}_{Tr}(q_0^1)$, which is well-defined by completeness of $Q$, compact, and it is independent of $n$. Therefore, the sequence $\{q_1^n\}$ is uniformly bounded over $[0, 1]$. Now observe that, for the Riemannian distance $d(\cdot, \cdot)$ and for all $0 \le t < \tau \le T$, $$d(q_1^n(t), q_1^n(\tau)) \le L\left(q_1^n\vert_{(t, \tau)}\right) = \int_t^{\tau}  \sqrt{g\left(\dot{q}_1^n(u), \dot{q}_1^n(u)\right)}du \le r(\tau - t),$$ where $q_1^n \vert_{(t,\tau)}$ denotes the restriction of $q_1^n$ to the interval $(t,\tau) \subset [0, T]$.
Therefore $\{q_1^n\}$ is equicontinuous on $[0, T]$, and by the Arzela-Ascoli Theorem, there then exists a subsequence $\{q_1^{n_k}\} \subset \{q_1^n\}$ which converges uniformly to a curve $q_1$ satisfying the boundary conditions in position. We may now replace $\{q^n\}$ with $\{q^{n_k}\}$ and repeat the analysis with $q_2^{n_k}$ to find some subsequence which converges uniformly to some $q_2$ satisfying the boundary conditions in position. Repeating this argument inductively over the agents, we hence obtain some subsequence of $\{q^n\}$—again denoted by $\{q^n\}$ for convenience—which converges uniformly to a curve $q = (q^1, q^2,..., q^s)$ satisfying the boundary conditions in position.

Let $(U_j^\mu, \varphi_j^\mu)$ be a finite collection of charts on $Q$ and $I_j^\mu$ an accompanying finite partition of $[0,1]$ such that, for sufficiently large $n$ and for all $1 \le j \le m$ and $\mu \in \mathcal{I}$, there exists a compact subset $K_j^\mu \subset U_j^\mu$ containing $q_j^n(I_j^\mu)$. In local coordinates, we may consider $q_j^n$ to be a curve on $\R^m$ (however, we will abuse this notation by continuing to call it $q_j^n$ both on the chart $U_j^{\mu}$ and its image in $\R^m$). Note that \begin{align}\label{christoffel}
\frac{d}{dt}\dot{q}_j^n = \frac{D}{\partial t} \dot{q}_j^n + \Gamma(q_j^n; \dot{q}_j^n, \dot{q}_j^n).
\end{align}where $\Gamma: \R^{3m} \to \R^m$ is continuous in the first argument and bilinear in the last two—and is determined by the ordinary Christoffel Symbols induced by the connection and chart. Hence, in $K_j^\mu$, we have 
\begin{align*}
    \left|\left|q_j^n \right|\right|_{H^2}^2 &= \left|\left| q_j^n \right|\right|_{L^2}^2 + \left|\left|\dot{q}_j^n \right|\right|_{L^2}^2 + \left|\left|\frac{d}{dt} \dot{q}_j^n \right|\right|_{L^2}^2 \\
    &= \int_{\mathcal{I}_j^\mu} \left|\left| q_j^n(t) \right|\right|_{\R^m}^2 dt + \int_{\mathcal{I}_j^\mu} \left|\left| \dot{q}_j^n(t) \right|\right|_{\R^m}^2 dt + \int_{\mathcal{I}_j^\mu}\left|\left|\frac{d}{dt} \dot{q}_j^n(t) \right|\right|_{\R^m}^2 dt.
\end{align*}

The first integral is bounded, as $\varphi_j^\mu$ is a homeomorphism and $K_j^\mu$ is compact—hence $\varphi_j^\mu(K_j^\mu) \supset (\varphi_j^\mu \circ q_j^n)(\mathcal{I}_j^\mu)$ is bounded in $\R^m$. For the second integral, note that for some scalars $\alpha, \ \beta$ we have $\alpha \left| \left| X \right| \right|^2_{\R^m} \le g(X, X) \le \beta \left| \left| X \right| \right|^2_{\R^m}$ for all $X \in T_x Q_j$ with $x \in K_j^\mu$. Hence, the boundedness of the second integral is equivalent to the boundedness of $$\int_{\mathcal{I}_j^\mu} g\left( \dot{q}_j^n(t), \dot{q}_j^n(t) \right)dt,$$ which follows by the uniform boundedness of $g(\dot{q}_j^n, \dot{q}_j^n)$ on $[0, 1]$ (and hence on the subset $\mathcal{I}_j^\mu$).

Similarly, the boundedness of the final integral is equivalent to the boundedness of  
$$\int_{I_j^{\mu}} g\left(\frac{d}{dt}\dot{q}_j^n(t), \frac{d}{dt}\dot{q}_j^n(t)\right)dt.$$

Note that $\Gamma$ and $V_{ij}$ are uniformly bounded on $I_j^{\mu}$ by continuity and the fact that each $q_j$ and $\dot{q}_j$ are uniformly bounded. It then follows by \eqref{christoffel} and by the fact that $\displaystyle{\sup_{n \in \N} } J(q^n) < +\infty$ that the above integral is bounded. Hence, $q_i^n$ is bounded in $H^2$. Since $H^2$ is a Hilbert space, we then get weak convergence of some subsequence of $(q^n)$ to $q \in \Omega$ in $H^2$.\hfill$\square$

\vspace{.2cm}

\begin{theorem}\label{PS}
$J$ satisfies the Palais-Smale condition on $\Omega$.
\end{theorem}
\textit{Proof:} Suppose that $\{q^n\} \subset \Omega$ is a Palais-Smale sequence for $J$, and let $B^n_i(t)$ and $C^n_i(t)$ be vector fields along $q_i^n$ satisfying $$\displaystyle{\frac{D^2}{\partial t^2} B^n_i(t) = R\left(\dot{q}_i^n(t), \frac{D}{\partial t}\dot{q}_i^n(t)\right)\dot{q}_i^n(t)} \ \hbox{ and  } \ \displaystyle{\frac{D^2}{\partial t^2} C^n_i(t) = \sum_{j \in \mathcal{N}_i} \grad_1 V_{ij}(q_i^n(t), q_j^n(t))}.$$ By Lemma \ref{Convergence}, we have that some subsequence of $q^n$ (again denoted by $q^n$) converges weakly to some curve $q\in \Omega$ with respect to the $H^2$ norm. Hence for all $i=1,...,s$, we have that $q^n_i$ and $\dot{q}^n_i$ are uniformly bounded, and $\frac{D}{\partial t}\dot{q}^n_i$ is bounded in the $L^2$ sense, so that $\displaystyle{R\left(\dot{q}_i^n(t), \frac{D}{\partial t}\dot{q}_i^n(t)\right)\dot{q}_i^n(t)}$ and $\displaystyle{\grad_1 V_{ij}(q_i^n, q_i^n)}$ are bounded in $L^2$. This implies, by using Gronwall's inequality, that $B^i_n$ and $C^i_n$ are bounded in $H^2$, so that there exists subsequences converging strongly in $L^2$.

Observe that $$\left|\left|\mathbf{d}J(q^n)\right|\right|_{T_{q^n}\Omega^\ast}\to 0 \implies \left|\left|\mathbf{d}_i J(q^n)\right|\right|_{T_{q_i^n}\Omega_i^\ast} =\left|\left|\grad_i J(q^n)\right|\right|_{T_{q_i^n}\Omega_i} \to 0,$$ for all $i = 1,...,s$, where $\mathbf{d}_i J$ represents the differential of the functional $J_i: \Omega_i \to \R$ defined by $J_i(p) = J(q_1,.., q_{i-1}, p, q_{i+1},...,q_s)$. Defining $A_i^n := \grad_i J(q^n)$, we have that $\frac{D^2}{\partial t^2}A^n_i \to 0$ in $L^2$, and for all $X \in T_{q_i^n}\Omega_i$,
\begin{align*}
    \mathbf{d}_i J(q^n)X &= \int_0^T g\left(\frac{D}{\partial t}\dot{q}_i^n(t) + B_i^n(t) + C_i^n(t), \frac{D^2}{\partial t^2} X\right)dt, \\
    \grad_i J(q^n)X &= \int_0^T g\left(\frac{D^2}{\partial t^2} A_i^n, \frac{D^2}{\partial t^2} X\right)dt.
\end{align*}
So that, with $\displaystyle{Z_i^n(t) := \frac{D}{\partial t}\dot{q}_i^n(t) + B_i^n(t) + C_i^n(t) - \frac{D^2}{\partial t^2} A_i^n}$, we have that for each $i = 1,...,s$ and $X \in T_{q_i^n}\Omega_i$,
\begin{align*}
    \int_0^T g\left(Z_i^n, \frac{D^2}{\partial t^2}X\right)dt = 0.
\end{align*}

Integrating by parts twice and applying the Fundamental Lemma of the Calculus of Variations, we see that $\frac{D^2}{\partial t^2} Z^n_i \equiv 0$. To conclude the result, it suffices to show that each $Z^n_i$ is bounded in $H^1$. This follows by defining the sequence of functions $\beta_i^n: \R \to \R$ given by $\beta_i^n(t) := g(Z_i^n(t), Z_i^n(t))$ and by noting that for all $t \in [0, T]$, $\beta_i^n(t) \ge 0$, $\frac{d^2}{dt^2}\beta_i^n \ge 0$ and $\frac{d^3}{dt^3}\beta_i^n \equiv 0$, from which it follows that $\beta_i^n(t)$ is just a sequence of quadratic polynomials on the real line. Moreover, since $\frac{d^2}{dt^2}\beta_i^n(t) = 2g\left(\frac{D}{\partial t} Z_i^n(t), \frac{D}{\partial t} Z_i^n(t)\right)$, the boundedness of $Z_i^n$ in the $H^1$ sense is equivalent to the boundedness of the sequence of coefficients to the quadratic term in $\beta_i^n(t)$. This follows from the fact that $\beta_i^n(t) \ge 0$ and $\frac{d^2}{dt^2}\beta_i^n \ge 0$ together with the boundedness of $\frac{D}{\partial t} Z_i^n$ on $L^2$. $\hfill\square$

As a consequence of Theorem \ref{PS} we have the following:

\begin{corollary}
There exists a curve $q = (q_1,...,q_s) \in \Omega$ such that each $q_i \in \Omega_i$ is smooth and satisfies \eqref{eqq1}, and such that $J(q) = \displaystyle{\inf_{\omega \in \Omega} J(\omega)}$.
\end{corollary}

\begin{remark}\label{potential}
In the case that the artifical potential function $V_{ij}$ is identically zero, the minimizers of $J$ are precisely the curves $q = (q_1,...,q_s) \in \Omega$ such that each $q_i$ is a Riemannian cubic polynomial satisfying the boundary conditions on $\Omega_i$ (see for instance \cite{marg}, \cite{CLACC}, \cite{CroSil:95}, \cite{Giambo}, \cite{noakes}).

One may wish to understand the influence of the potential in shaping the minimizers. To this end, observe that if $q_i(t)$ is a Riemannian cubic and $\displaystyle{X_i = -\sum_{j \in \mathcal{N}_{i}}\grad_1 V_{ij}(q_i, q_j)}$, for $i = 1,...,s$, we have \begin{align*}
    \mathbf{d}_i J(q)X_i = -\int_0^T g(X_i, X_i) dt \le 0.
\end{align*} Hence, $J$ decreases when one deforms the Riemannian cubic $q_i$ in the direction that the artificial potential is decreasing fastest. 

For the application of collision avoidance, we desire that the minimizers are given by deforming the cubics in a way that the average distance between agents is increased. In other words, we would like $J$ to be decreasing along the sum of the gradient flows of the Riemannian distances between the agents:

$$\mathbf{d}_i J(q) \sum_{(i, j) \in \mathcal{E}} \grad_1 d(q_i, q_j) = \int_0^T g\left(\sum_{(i, j) \in \mathcal{E}}\grad_1 V_{ij}(q_i, q_j), \sum_{(i, j) \in \mathcal{E}}\grad_1 d(q_i, q_j)\right)dt \le 0.$$

\vspace{.1cm}

Note that any potential whose gradient is of the form $\sum_{(i,j)\in\mathcal{E}}\grad_1 V_{ij}(p, q)=$\newline  $-f(p, q)\sum_{(i,j)\in\mathcal{E}}\grad_1 d(p, q)$ will work, where $f: Q \times Q \to \R$ is non-negative. A particularly simple, smooth, positive-definite family of such potentials is given by $\displaystyle{V_{ij}= \frac{1}{\epsilon + (d(p, q)/D)^k}}$ for $\epsilon, D > 0,$ and $k \in \mathbb{N}$ and for all $(i, j) \in \mathcal{E}$.\end{remark}

\section{The Collision Avoidance Task}\label{Sec: Collision Avoidance}
In this section, we focus on the task of collision avoidance between nearest neighbors. First we define precisely what it means for collision to be avoided, and prove that the minimizers of $J$ avoid collision for a general class of artificial collision avoidance potentials. The strategy proposed is to use a reference trajectory $q \in \Omega$ that avoids collision, from which we obtain an upper bound on $J(q)$—which must also be an upper bound for the value of $J$ at a minimizing trajectory.

\subsection{Collision avoidance}\label{avoidance}
For $(i, j) \in \mathcal{E}$ and real numbers $0 < r_{ij} < r_{ij}^\ast < R_{ij}$, we define the sets:
\begin{align*}
    C_{ij} &:= \{(p, q) \in Q \times Q \ : \ d(p, q) < r_{ij} \} \quad &&\textit{Collision Region } \text{for the edge 
    $(i,j)\in\mathcal{E}$,} \\
    C^\ast_{ij} &:= \{(p, q) \in Q \times Q \ : \ d(p, q) < r_{ij}^\ast \} &&\textit{Risk Region } \text{for the edge
    $(i,j)\in\mathcal{E}$,}\\
    S_{ij} &:= \{(p, q) \in Q \times Q \ : \ d(p, q) > R_{ij} \} &&\textit{Safety Region } \text{for the edge 
    $(i,j)\in\mathcal{E}$,}
\end{align*}where $r_{ij}, \ r_{ij}^\ast,\  R_{ij}$ are called the \textit{tolerances} for the collision region, risk region, and safety region, respectively. We assume that each of these regions is symmetric in $(i, j)$. That is, $C_{ij} = C_{ji}$ and equivalent statements for $C_{ij^\ast}$ and $S_{ij}$. We say that $q \in \Omega$ \textit{avoids collision with tolerances $r_{ij}$} if $(q_i(t), q_j(t)) \notin C_{ij}$ for all $(i. j) \in \mathcal{E}$ and $t \in [0, T]$. In Figure \ref{fig: avoid_reg} below, we depict these regions with respect to agent $i$.

\begin{figure}[h!]\label{fig: avoid_reg}
\begin{center}
 \includegraphics[width=6.cm]{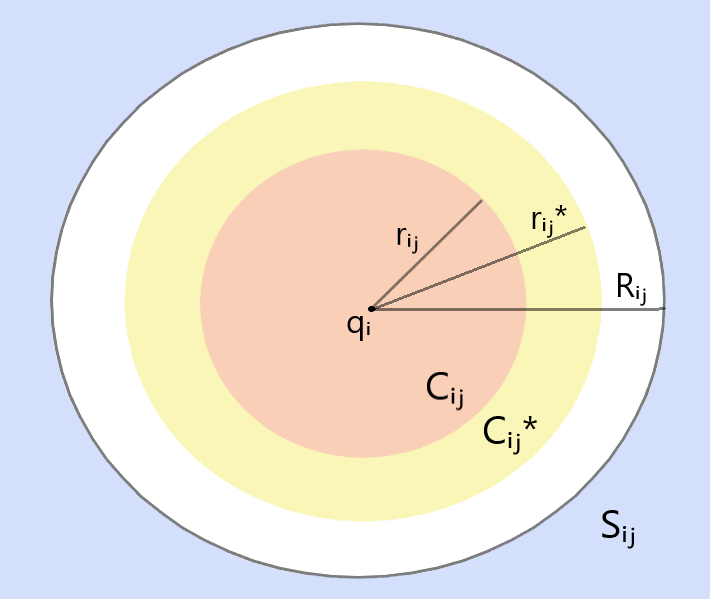}
 \caption{The regions $C_{ij}, C_{ij}^\ast$, and $S_{ij}$ taken with respect to the position of agent $i$. For example, collision between agents $(i, j) \in \mathcal{E}$ is equivalent to saying that the position of agent $j$ enters the red ball somewhere along the trajectories.}
 \label{figuav}
 \end{center}
\end{figure}

We will construct our potential functions for collision avoidance so that its components are bounded above in the corresponding Safety Regions and bounded below in the corresponding Risk Region.  More precisely, for all $(i, j) \in \mathcal{E}$ and for some real numbers $0 < V_{ij}^- < V_{ij}^\ast$, we construct the artificial potential $V_{ij}$ such that $V_{ij} > V_{ij}^\ast$ on $C_{ij}^\ast$ and $V_{ij} < V_{ij}^{-}$ on $S_{ij}$.
From equation \eqref{eqq1}, minimizers of $J$ will move in a way that the fourth covariant derivatives of the agents follow the gradient descent of the potential (perturbed by some quantity that depends on the geometry of the manifold). One would thus expect the minimizers to avoid collision if $V_{ij}^\ast - V_{ij}^-$ is sufficiently large. Indeed this is the case (with some caveats dependent on the geometry of $\Omega$ and the chosen tolerances), but in order to more precisely define "sufficiently large," one must note the influence of the boundary conditions and the tolerances of the Collision regions, Risk regions, and Safety regions. In particular, the rate at which the agents must accelerate away from the Collision region once they have entered the Risk region will depend on the quantity $r_{ij}^\ast - r_{ij}$ together with the speeds at which the agents enter the Risk region. In turn, these speeds will be determined indirectly by equation \eqref{eqq1} together with the boundary conditions, which depend upon the magnitude of the potential's gradient along the trajectories (which itself may be influenced by $V_{ij}^\ast - V_{ij}^-$). This circular dependence poses a problem in the analysis that we resolve by introduce some "extra data" that can be indirectly used to bound the troublesome quantities. 

Consider a \textit{reference trajectory} $q = (q_1,...,q_s) \in \Omega$ such that $(q_i(t), q_j(t)) \in S_{ij}$ for each $(i, j) \in \mathcal{E}$ and $t \in T$ (note that this requires $R_{ij} < \min\left\{d(q_i^0, q_j^0), \ d(q_i^T, q_j^T)\right\}$). We define for all $i\in\mathcal{V}$ the non-negative real numbers:
\begin{equation*}
    a_i := \sup_{t \in [0, T]} \left\{ \left|\left|\frac{D}{\partial t}\dot{q}_i(t)\right|\right|\right\},\quad
    c := T\sum_{i=1}^s \left[a_i^2 + \sum_{j \in \mathcal{N}_i}V_{ij}^{-} \right],\quad
    v_i := \sqrt{cT} + \sqrt{cT + ||v_i^0||^2}.
\end{equation*}

\begin{proposition}\label{Collision_Avoidance_Prop}
If $V_{ij}^\ast > \frac{c(v_i + v_j)}{2(r_{ij}^\ast - r_{ij})}$ for all $(i,j)\in\mathcal{E}$, then every minimizer $q^\ast = (q_1^\ast,...,q_s^\ast) \in \Omega$ of $J$ avoids collision.
\end{proposition}

\textit{Proof:} Observe that, since each pair of agents in our reference trajectory remains in the Safety Region for all $t \in [0, T]$,
\begin{align*}
    J(q) &= \sum_{i=1}^s\int_0^T\left[ \left|\left|\frac{D}{\partial t}\dot{q}_i(t) \right|\right|^2 + \sum_{j \in \mathcal{N}_i}V_{ij}(q_i(t), q_j(t)) \right]dt \\
    &\le T\sum_{i=1}^s \sup_{t \in [0, T]}\left[ \left|\left|\frac{D}{\partial t}\dot{q}_i(t) \right|\right|^2 + \sum_{j \in \mathcal{N}_i}V_{ij}(q_i(t), q_j(t)) \right] \\
    &\le T\sum_{i=1}^s \left[a_i^2 + \sum_{j \in \mathcal{N}_i}V_{ij}^{-} \right] =: c
\end{align*}
For all $i \in \mathcal{V}$, let $v_i^\ast = \displaystyle{\sup_{t \in [0, T]} \left\{||\dot{q}^\ast_i(t)|| \right\}}$, so that
\begin{align*}
   ||\dot{q}^\ast_i(t)||^2 &= ||v_i^0||^2 + \int_0^T \frac{d}{dt}||\dot{q}^\ast_i(t)||^2 dt \\
   &= ||v_i^0||^2 + 2\int_0^T g\left(\frac{D}{\partial t}\dot{q}^\ast_i(t), \dot{q}^\ast_i(t)\right)dt \\
   &\le ||v_i^0||^2 + 2\left[\int_0^T \left|\left|\frac{D}{\partial t}\dot{q}^\ast_i(t) \right|\right|^2 dt\right]^{1/2} \left[\int_0^T \left|\left|\dot{q}^\ast_i(t) \right|\right|^2 dt \right]^{1/2} \\
   &\le ||v_i^0||^2 + 2v_i^\ast\sqrt{J(q^\ast)T} \\
   &\le ||v_i^0||^2 + 2v_i^\ast\sqrt{cT},
\end{align*}
since $J(q^\ast) \le J(q) \le c$. Taking the supremum over $t \in [0, T]$ on both sides of the inequality, we have $(v_i^\ast)^2 \le ||v_i^0||^2 + 2v_i^\ast\sqrt{cT}$, so that $v_i^\ast \le \sqrt{cT} + \sqrt{cT + ||v_i^0||^2} =: v_i$.

Now assume towards contradiction that $q^\ast$ does not avoid collision. Then for some $(i,j) \in \mathcal{E}$ and $\tau \in [0, T]$, $(q_i(\tau), q_j(\tau)) \in C_{ij}$, and so for $t > \tau$ we have by repeated applications of the triangular inequality that

\begin{align*}
    d(q_i^\ast(t), q_j^\ast(t)) &\le d(q_i^\ast(t), q_i^\ast(\tau)) + d(q_i^\ast(\tau), q_j^\ast(t)) \\
    &\le d(q_i^\ast(t), q_i^\ast(\tau)) + d(q_j^\ast(\tau), q_j^\ast(t)) + d(q_i^\ast(\tau), q_j^\ast(\tau)) \\
    &\le L\left(q_i|_{(\tau, t)} \right) + L\left(q_j|_{(\tau, t)} \right) + r_{ij} \\
    &= \int_{\tau}^t ||\dot{q}_i(t)||dt + \int_{\tau}^t ||\dot{q}_j(t)||dt + r_{ij} \\
    &\le (v_i^\ast + v_j^\ast)(t - \tau) + r_{ij} \\
    &\le (v_i + v_j)(t - \tau) + r_{ij},
\end{align*}
and so for all $t \in \left[\tau, \tau + \frac{r_{ij}^\ast - r_{ij}}{v_i + v_j}\right)$ then \ $((q_i(t), q_j(t)) \in C_{ij}^\ast$. It can be seen by a similar argument that $((q_i(t), q_j(t)) \in C_{ij}^\ast$ for all $t \in \left(\tau - \frac{r_{ij}^\ast - r_{ij}}{v_i + v_j}, \tau\right]$. Therefore,

\begin{align*}
    2V_{ij}^\ast \left(\frac{r_{ij}^\ast - r_{ij}}{v_i + v_j} \right) &\le \int_{\tau - \frac{r_{ij}^\ast - r_{ij}}{v_i + v_j}}^{\tau + \frac{r_{ij}^\ast - r_{ij}}{v_i + v_j}} \sum_{j \in \mathcal{N}_i} V_{ij}(q_i(t), q_j(t))dt\\
    &\le \int_0^T \sum_{j \in \mathcal{N}_i} V_{ij}(q_i(t), q_j(t))dt \\
    &\le J(q^\ast) \le J(q) \le c.
\end{align*}

So that $V_{ij}^\ast \le  \frac{c(v_i + v_j)}{2(r_{ij}^\ast - r_{ij})}$. Therefore, by contradiction, $q^\ast$ avoids collision. \hfill$\square$


In the context of Proposition \ref{Collision_Avoidance_Prop}, one may think of the reference trajectory $q$ as avoiding collision with tolerances $R_{ij}$. To that end, we call collision avoidance with tolerances $r_{ij}$ \textit{feasible} if there exists a curve $q \in \Omega$ which avoids collision with tolerances $r_{ij}$.

We now consider the smooth family of repulsive potentials parameterized by $D,\epsilon \in \R^+, \ k \in \N$ defined by \begin{equation}\label{Vtune}V_{D, \epsilon}^k(p, q) = \frac1{\epsilon + (d(p, q)/D)^k}.\end{equation}

\begin{corollary}\label{cor: Coll_Avoid}
For all $(i,j) \in \mathcal{E}$, if collision avoidance is feasible for the tolerances $R_{ij} > 0$, then for all $r_{ij} < R_{ij}$, there exists $\epsilon_{ij} \in \R^+, \ k_{ij} \in \N$ such that every minimizer of $J$ with the potential $V_{ij} = V_{R_{ij}, \epsilon_{ij}}^{k_{ij}}$ avoids collision with tolerances $r_{ij}$. 
\end{corollary}

\textit{Proof:}
Let $r_{ij} < R_{ij}$, and choose $r_{ij}^\ast$ such that $r_{ij} < r_{ij}^\ast < R_{ij}$. Since collision avoidance is feasible with the tolerances $R_{ij}$, there exists a reference trajectory $q \in \Omega$ such that each pair of agents in $\mathcal{E}$ remains in the safety region with tolerances $R_{ij}$.

We have $V_{ij}(p, q) < 1$ whenever $d(p, q) \ge R_{ij}$, independent of $\epsilon$ and $k$, so that we may suppose $V_{ij}^- = 1$. Moreover, $d(p, q)/R_{ij} < r_{ij}^\ast/R_{ij} < 1$ whenever $d(p, q) < r_{ij}^\ast$, so that $(d(p, q)/D_{ij})^{k_{ij}}$ can be made arbitrarily small by taking $k_{ij}$ sufficiently large—in which case $V_{ij}(p, q)$ becomes arbitrarily close to $\frac1{\epsilon_{ij}}$. In particular, we may suppose that $V_{ij}^\ast = \frac{1}{2\epsilon_{ij}} < \frac1{\epsilon_{ij}}$. Finally, observe that $\frac{c(v_i + v_j)}{2(r_{ij}^\ast - r_{ij})}$ is finite and independent of the parameters $\epsilon_{ij}, k_{ij}, D_{ij}$. Hence, by choosing $\epsilon_{ij} < \frac{r_{ij}^\ast - r_{ij}}{c(v_i + v_j)}$ we have by Proposition \ref{Collision_Avoidance_Prop} that any minimizer of $J$ avoids collision with tolerances $r_{ij}$. \hfill$\square$

\begin{remark}
Notice that the family of potentials defined above is strictly positive, so that we are assuming that the neighboring agents can always measure their distances with respect to one another. In practice, this is often not the case due to technological limitations. Instead, we have some \textit{sensing radius} which dictates how close two agents must be to sense one another (i.e., to measure their distance). When two agents are outside of the sensing radius, the component of the artificial potential relating these agents must not include their distance, and in the simplest case is constant.

In order to preserve the regularity of the potential (recall that $C^1$ was required in the proof of Theorem \ref{PS} for the existence of minimizers), we can define a bump function. For example, consider
$$V(p, q) = \begin{cases}
\frac{1}{\epsilon}\exp\left(-\frac{1}{1 - (d(p, q)/D)^k}\right) & d(p, q) < D \\
0 & \text{else}.
\end{cases}$$

$V$ is smooth and can easily replace the family of potentials used in Corollary \ref{cor: Coll_Avoid} to achieve the same result, except with collision avoidance being guaranteed for all $r_{ij} < \min\{R_{ij}, h\}$.
\end{remark}

\subsection{Collision avoidance with uniformly bounded derivatives}\label{Sec: Avoidance_bounded_domain}

In the path planning design of some robotic applications sometimes it is needed to require uniform bounds on the magnitudes of the derivatives of the trajectories (see \cite{jon} and  \cite{mac} for instance). While the existence of minimizers is not guaranteed in such a case, we still find it instructive to provide analysis in such a case. In particular, we provide alternate conditions for $V_{ij}^\ast$—expressed in terms of the uniform bounds and initial conditions—which ensures the collision avoidance of minimizers. We do so by taking advantage of an invariant for the flow of equations \eqref{eqq1}, and hence the analysis will show collision avoidance for all critical points of $J$ rather than just the minimizing trajectories (as was the case in Section \ref{avoidance})—provided that the derivatives of these solutions satisfy the given bounds. Moreover, as will be seen below, the conditions on the potential no longer have any dependence on the Risk and Safety regions, so that we may suppose that $r_{ij}^\ast = r_{ij}$ and set $V_{ij}^- = V_{ij}(q_i^0, q_j^0)$.

For each $i \in \mathcal{V}$, we define the domain $$\Omega_i^{\max} = \left\{q \in \Omega_i \ | \ \left|\left|\dot{q}(t)\right|\right| \le v_i^{\max}, \ \left|\left|\frac{D}{\partial t} \dot{q}(t)\right|\right| \le a_i^{\max}, \ \left|\left|\frac{D^2}{\partial t^2} \dot{q}(t)\right|\right| \le \eta_i^{\max}, \ \  \forall t \in [0, T] \right\},$$ 
where $v_i^{\max}, a_i^{\max}, \eta_i^{\max}$ are non-negative real numbers with $v_i^{\max} > \max \left\{\left|\left|v_i^0\right|\right|, \ \left|\left|v_i^T\right|\right| \right\}$ and further define $\Omega^{\max} := \Omega_1^{\max} \times \cdots \times \Omega_s^{\max}$. We further set $\eta_i^0 = \frac{D^2}{\partial t^2} \dot{q}(0)$.

\begin{proposition}\label{Prop: Coll_avoid_bounded}
Suppose that $q = (q_1,...,q_s) \in \Omega^{\max}$ and that each $q_i$ satisfies equation $\eqref{eqq1}$. If $$V^{\ast}_{ij} > \sum_{i=1}^s \left[(a_i^{\max})^2 + v_i^{\max}\eta_i^{\max} + v_i^0 \eta_i^0 + \frac12 \sum_{j \in \mathcal{N}_i} V_{ij}^{-} \right]$$ for all $(i, j) \in \mathcal{E}$, then $q$ avoids collision with tolerances $r_{ij}$.
\end{proposition}

\textit{Proof:} We define the function $H$ by:
$$H\left(q, \dot{q}, \frac{D}{\partial t} \dot{q}, \frac{D^2}{\partial t^2} \dot{q}\right) = \sum_{i=1}^s\left[
g\left(\dot{q}_i, \frac{D^2}{\partial t^2}\dot{q}_i\right) 
- \frac12 \left|\left|\frac{D}{\partial t} \dot{q}_i\right|\right|^2
+ \frac12\sum_{j \in \mathcal{N}_i} V_{ij}(q_i, q_j) \right].$$
Using the symmetry of the potential and since $\mathbb{G}$ is undirected, we have that the time evolution of $H$ satisfies
\begin{align*}
    \dot{H} &= \sum_{i=1}^s \Bigg[ - g\left(\dot{q}_i, R\left(\dot{q}_i, \frac{D}{\partial t} \dot{q}_i\right)\dot{q}_i + \sum_{j \in \mathcal{N}_i} \grad_1 V_{ij}(q_i, q_j) \right)\\& \hspace{1.5cm}+ \frac12\sum_{j \in \mathcal{N}_i} \left(g\left(\grad_1 V_{ij}(q_i, q_j), \dot{q}_i\right) + g\left(\grad_j V_{ij}(q_i, q_j), \dot{q}_j\right)\right) \Bigg] \\
    &= \sum_{i=1}^s \left[-g\left(\dot{q}_i, R\left(\dot{q}_i, \frac{D}{\partial t} \dot{q}_i\right) \dot{q}_i  \right) - \sum_{j \in \mathcal{N}_i} g\left(\dot{q}_i, \grad_1 V_{ij}(q_i, q_j)\right) + \sum_{j \in \mathcal{N}_i} g\left(\grad_1 V_{ij}(q_i, q_j), \dot{q}_j\right) \right] \\
    &= -\sum_{i=1}^s Rm\left(\dot{q}_i, \frac{D}{\partial t} \dot{q}_i, \dot{q}_i, \dot{q}_i\right),
\end{align*}
which vanishes identically due to the symmetries of the curvature tensor (see Theorem 3.1 in \cite{Boothby} for more details). Suppose that $q$ does not avoid collision. Then, there exists $(i, j) \in \mathcal{E}$ and some $\tau \in [0, T]$ such that $d(q_i(\tau), q_j(\tau)) \le r_{ij}$, so that $V_{ij}(q_i(\tau), q_j(\tau)) > V_{ij}^{\ast}$. But, since $\dot{H} \equiv 0$, we have:
\begin{align*}
    V^{\ast}_{ij} - \sum_{i=1}^s\left[\frac12 (a_i^{\max})^2 + v_i^{\max}\eta_i^{\max}\right] &\le H(\tau) = H(0) \le \sum_{i=1}^s \left[v_i^0 \eta_i^0 + \frac12 \sum_{j \in \mathcal{N}_i} V_{ij}^{-} \right] \\
    \implies V^{\ast}_{ij} &\le \sum_{i=1}^s \left[(a_i^{\max})^2 + v_i^{\max}\eta_i^{\max} + v_i^0 \eta_i^0 + \frac12 \sum_{j \in \mathcal{N}_i} V_{ij}^{-} \right] 
\end{align*}
Which contradicts the assumption on $V_{ij}^{\ast}$. $\hfill\square$

\subsection{Simulation results}
In this section, we conduct numerical simulations of the collision avoidance problem in the cases of $Q = \R^3$ and $Q = S^2$. In both cases, numerical integration was done via the Euler method with a time step of $h = 0.005$. Initial conditions were decided so that the unique geodesics satisfying those conditions collide along the trajectories. The final points of the geodesics were then used as the end points in the boundary value problems satisfying \eqref{eqq1} along the trajectories. A shooting method based on the downhill simplex algorithm was used to find the initial accelerations and jerks that lead to solutions to the boundary value problem. 
\subsubsection{Collision Avoidance on $\R^3$}
Consider 4 agents on $\R^3$ with the Euclidean metric, and the following neighboring relations between agents $\mathcal{N}_1 = \{2, 3, 4\}, \ \mathcal{N}_2 = \{1, 3\}, \ \mathcal{N}_3 = \{1,2,4\}, \ \mathcal{N}_4 = \{1, 3\}$. We choose $T = 4$ and the boundary conditions as:
\begin{align*}
    &q_1(0) = (-1, -1, 0), &&q_1(T) = (1, 1, 0), &&\dot{q}_1(0) = \dot{q}_1(T) = (0.5, 0.5, 0)\nonumber \\
    &q_2(0) = (-1, 1, 0), &&q_2(T) = (1, -1, 0), &&\dot{q}_2(0) = \dot{q}_2(T) = (0.5, -0.5, 0)\nonumber \\
    &q_3(0) = (1, 1, 0), &&q_3(T) = (-1, -1, 0), &&\dot{q}_3(0) = \dot{q}_3(T) = (0.5, 0.5, 0) \nonumber\\
    &q_4(0) = (1, -1, 0), &&q_4(T) = (-1, 1, 0), &&\dot{q}_4(0) = \dot{q}_4(T) = (-0.5, 0.5, 0)\nonumber \\
\end{align*}

For our reference trajectories, we choose piecewise defined curves consisting of cubic polynomials, with constant speed arcs along the circle of radius $2$ centered at the origin. For example, the reference trajectory $\gamma_1$ for $q_1$ consists of the cubic polynomial verifying the relevant initial conditions and $\gamma(1.5) = (2, 0, 0), \ \dot{\gamma}(1.5) = (0, -\pi, 0)$, followed by the circular arc $\gamma(t) = (2\cos(\frac{\pi}{2}(t - 1.5)), \ -2\sin(\frac{\pi}{2}(t - 1.5),\ 0)$ for $1.5 < t < 2.5$ , followed by the cubic polynomial starting from the circular arcs end points in position and velocity to the final points in the boundary value problem. Note that, in this manner $\gamma_1 \in \Omega_1$. The other reference trajectories were chosen similarly, so that $\gamma = (\gamma_1,..., \gamma_4) \in \Omega$. It can be seen that $||\gamma_i(t) - \gamma_j(t)|| > 2$ for all $(i, j) \in \mathcal{E}$ and $t \in [0, 4]$. Consistent with Section \ref{avoidance}, we choose the repulsive potential with $V_{ij}(p, q) = \frac{1}{\epsilon + (||p - q||/2)^k}$ so that $V_{ij}^- = 1$.

Moreover, it can be shown numerically that $a_i < 3.22$, so that $c < 185.9$ and $v_i < 54.57$ for $i = 1,...,4$. We choose our tolerances for the Collision regions and Risk regions as $r_{ij} = \frac12$ and $r_{ij}^\ast = 1$ for all $(i, j) \in \mathcal{E}$, so that by Proposition \ref{Collision_Avoidance_Prop}, collision is avoided if we choose $k$ and $\epsilon$ such that $\frac1{\epsilon + 2^{-k}} > 20,290$. One such solution is $\epsilon = 0.00001$ and $k = 16$, with which the necessary conditions for extrema \eqref{eqq1} take the form

\begin{align}\label{Ex: R3 nec}
    \frac{d^4 q_i}{dt^4}  = \sum_{j \in \mathcal{N}_i} \frac{64||q_i - q_j||^{14}}{(0.00001 + ||q_i - q_j||^{16})^2}(q_i - q_j), \qquad i=1,\ldots,4.
\end{align}

In Figure \ref{fig: R3} below, we show the numerical solution to equations \eqref{Ex: R3 nec} with the above initial conditions. For the purposes of comparison, we also include the numerical solutions to the cubic polynomials satisfying the boundary conditions (that is, the solutions to the necessary conditions for extrema in the case that $V \equiv 0$).

\begin{figure}[h!]
\begin{center}
 \includegraphics[width=7.cm]{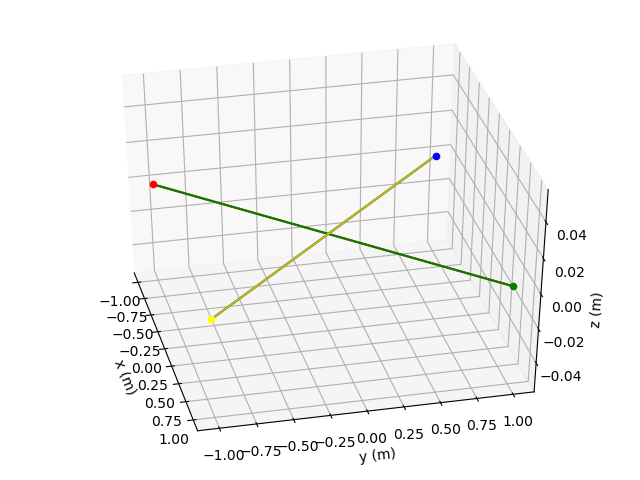}
 \includegraphics[width=6.6cm]{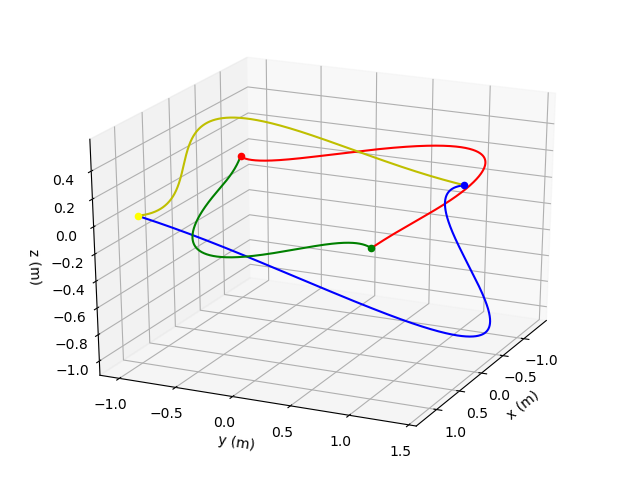}
 \caption{The Riemannian cubics (left) and solutions to the necessary conditions for extrema (right) satisfying the boundary conditions in positions and velocities as above. In the case of Riemannian cubics, all four agents collide at the center. With the potential, it is easily verified that collision is avoided within the desired tolerance given by $r_{ij} = \frac12$.}\label{fig: R3}
 \end{center}
\end{figure}



\subsubsection{Collision Avoidance on $S^2$}

Denote by $H:=G/K$ a Riemannian symmetric space, where $G$ is a compact and connected  finite-dimensional Lie group endowed with a bi-invariant Riemannian metric
and   $K$ a closed Lie subgroup of $G$. The canonical projection $\pi:G\to H$ is a Riemannian submersion (see \cite{Helgason} for instance). Therefore, for all $g$ in $G$, the isomorphism $T_g\pi:(\ker T_g\pi)^{\bot} \to T_{\pi(g)}H$ preserves the inner-products defined by the Riemannian metrics on $G$ and $H$, and $T_gG$ splits into, the vertical subspace $V_g:=\ker T_g\pi$ and the horizontal subspace $\hbox{Hor}_g=(V_g)^{\perp}:=(\ker T_g\pi)^{\bot}$. 
In particular,  the Lie algebra $\mathfrak{g}$ of $G$ admits the decomposition
$\mathfrak{g}=\mathfrak{s}\oplus\mathfrak{m}$ where $\mathfrak{s}$ is the Lie algebra of $K$ and $\mathfrak{m}\simeq T_{\mathfrak{o}}H$, with  $\mathfrak{o}=\pi(e)$, and $e$ the identity element on $G$. That is, $\ker T_e\pi=\mathfrak{s}$  and the horizontal subspace $(\ker T_g\pi)^{\bot}$ is $\mathfrak{m}$. Moreover, the relations
$[\mathfrak{s},\mathfrak{s}]\subset\mathfrak{s},\,[\mathfrak{m},\mathfrak{m}]\subset\mathfrak{s},\,[\mathfrak{m},\mathfrak{s}]\subset\mathfrak{m}$, hold (see \cite{Helgason}, Ch. IV, Sec. 5). Using this decomposition of $T_gG$, it is possible to extend the notion of vertical and horizontal tangent vectors on $G$ to vertical and horizontal vector fields and curves (see for instance \cite{ZN}).

We return for a moment to the family of potential functions  considered in Section \ref{avoidance} given by $V_{\epsilon, D}^k(p, q) = \frac{1}{\epsilon + (d(p,q)/D)^k}$ (see equation \eqref{Vtune}). Recall that when $p, q$ are sufficiently close, the Riemannian distance function can be written in terms of the Riemannian exponential $\exp_q$ as $d(p, q) = \left|\left| \exp_p^{-1}q\right|\right|$. This representation provides us with a way to calculate the gradient of the potential, which is summarized in the following lemma.

\begin{lemma}
Suppose that $Q$ is a Riemannian manifold and $p, q$ are contained within some convex ball on which $\exp_p$ is a diffeomorphism. Then $$\grad_1 V_{\epsilon, D}^{k}(p, q) = \frac{k d(p,q)^{k-2}}{2D^k (\epsilon + (d(p,q)/D)^k)^2} \exp_p^{-1} q.$$
\end{lemma}

\textit{Proof:} It is shown in \cite{point} that $\grad_1 d^2(p, q) = \exp_p^{-1}q$. The result follows immediately by writing $d(p, q)^k = d^2(p, q)^{k/2}$ and applying the chain rule.\hfill$\square$



In particular, when $Q = G$ is a Lie group, the Riemannian exponential at the point $g \in G$ can be represented by the Riemannian exponential at the identity element $e \in G$ as $\exp_g^{-1}(h) = \exp_e^{-1}(g^{-1}h)$. Now, similar to \cite{point} and \cite{CoGo20}, if we consider the horizontal curve $\tilde{g}_i\in\pi\mid_{G}^{-1}(g_j)$ on $G$, $i\in\mathcal{V}$, verifying  $\displaystyle{\xi_i=T_{\tilde{g}_i(t)}L_{\tilde{g}_i^{-1}}\left(\frac{d\tilde{g}_i}{dt}\right)}$ and
\begin{equation}\label{eq1withexp4}
\xi'''_i=-[\xi_i, [\xi^{\prime}_i,\xi_i]]-\frac{1}{2}\sum_{j\in\mathcal{N}_i}\frac{k_{ij}||\exp_e^{-1}(\tilde{g}^{-1}_i\tilde{g}_j)||^{k_{ij}-2}}{D_{ij}^k(\epsilon_{ij} + (||\exp_e^{-1}(\tilde{g}^{-1}_i\tilde{g}_j)||)/D)^{k_{ij}})^2}\exp_e^{-1}(\tilde{g}^{-1}_i\tilde{g}_j),\quad i\in\mathcal{V},
\end{equation}
where the latter equation evolving on the subspace $\mathfrak{m}$, then the solutions project down to the solutions to the necessary conditions \eqref{eqq1} on $H$ with the potential defined by $V_{ij} = V_{\epsilon_{ij}, D_{ij}}^{k_{ij}}$.

For simulation purposes we now restrict to the case of three identical agents on the symmetric space $H = S^2$, the two-dimensional unit sphere, with $G = \SO(3)$ and $K = \SO(2)$, for $i = 1, 2, 3$. Denoting the canonical basis of $\R^3$ by $\{e_1, e_2, e_3\}$, the Lie group $SO(2)$ can be seen as the subgroup of $\SO(3)$ leaving $e_1 \in S^2$ fixed. In such a case, the Riemannian exponential map $\text{exp} : \so(3) \to \SO(3)$ is just the usual matrix exponential map
on $SO(3)$, where $\mathfrak{so}(3)$ is the Lie algebra of $SO(3)$, i.e., the set of all $3\times3$ skew-symmetric matrices. 

The matrix exponential map is a diffeomorphism between $\mathcal{U} = \{\hat{a} \in \so(3): a \in \R^3, ||a|| < \pi\}$ and $\mathcal{K} = \{R \in \SO(3) : \text{Tr}(R)\ne -1\}$, and its inverse map is the matrix logarithm map.
Each agent $q_i$ on $S^2$ can be represented by $R_i$ on $\SO(3)$ via the relation $q_i = R_i e_1$ and the projection $\pi: \SO(3) \to S^2$ given by $\pi(R_i) = R_i e_1$. Denote by $\phi=\hbox{arccos}(\frac{1}{2}(\hbox{Tr}(R_i^{T}R_j)-1)$, and using Proposition $5.7$ in \cite{bookBullo},  for $R_i\neq R_j$ then $\log(R_i^{T}R_j)=\frac{\phi}{2\sin(\phi)}(R_i^{T}R_j-R_j^{T}R_i)$ and $\|\log(R^{T}_iR_j)\|=\phi$. For all $(i, j) \in \mathcal{E}$, we choose the parameters $k_{ij} = 8, \ D_{ij} = \frac12, \ \epsilon_{ij} = 0.00001$ for our potential. If we denote by $\times:\mathfrak{so}(3)\to\mathbb{R}^{3}$ the inverse of the hat isomorphism $\hat{\cdot}:\mathbb{R}^3\to\mathfrak{so}(3)$, it follows that the necessary conditions for extrema \eqref{eq1withexp4} are given by:

\begin{equation}\label{eqwithexp4}
\xi_i'''=-\xi_i \times (\xi_i^{\prime}\times \xi_i)-\sum_{j\in\mathcal{N}_i}\frac{512\phi(R_i^T R_j)^{7}(R_iR_j^{T}-R_jR_i^{T})^{\times}}{\sin(\phi(R_i^T R_j))( 0.00001 + 256\phi(R_i^T R_j)^8)^2}
\end{equation} together with the equation $\dot{R}_i=R_i\hat{\xi}_i$, the condition $\pi(R_j)=q_j$, and the boundary conditions $R_i(0)=R_0^i$, $R_i(T)=R_T^{i}$, $\xi_i(0)=\xi_0^i$, $v_i(T)=v_T^{i}$. We then obtain the solution $q_i$ to the corresponding necessary conditions on $H$ by projecting $R_i$ to $S^2$.

We set $T = 4$, and chose initial conditions as:
\begin{align*}\label{Ex: S2 BCs}
    &R_1(0) = \begin{bmatrix} 1 & 0 & 0 \\ 0 & 1 & 0 \\ 0 & 0 & 1\end{bmatrix}, &&R_2(0) = \begin{bmatrix} 0 & -1 & 0 \\ 1 & 0 & 0 \\ 0 & 0 & 1\end{bmatrix}, && R_3(0) = \begin{bmatrix} \frac{\sqrt{2}}{2} & -\frac{\sqrt{2}}{2} & 0 \\ \frac{\sqrt{2}}{2} & \frac{\sqrt{2}}{2} & 0 \\ 0 & 0 & 1 \end{bmatrix} \\
    &\xi_1(0) = [0, 0.24, 0.25]^T, && \xi_2(0) = [0, 0.2, -0.25]^T, &&\xi_3(0) = [0, 0.2, 0.02]^T 
\end{align*}

Figure \ref{fig: s2} below shows the geodesics with these initial conditions and the solutions to the necessary conditions \eqref{eqwithexp4} satisfying the same boundary conditions as the geodesics.

\begin{figure}[h!]
\begin{center}
 \includegraphics[width=6.cm]{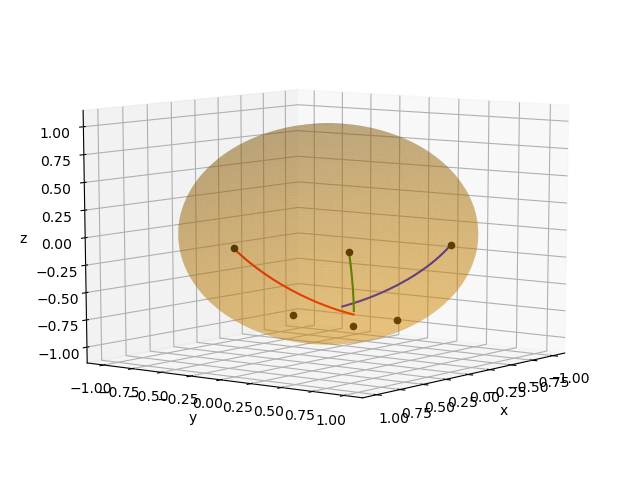}
 \includegraphics[width=6.cm]{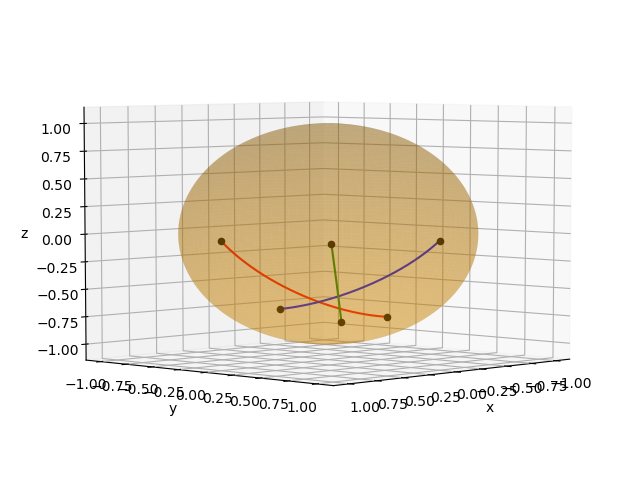} \\
 \includegraphics[width=6.cm]{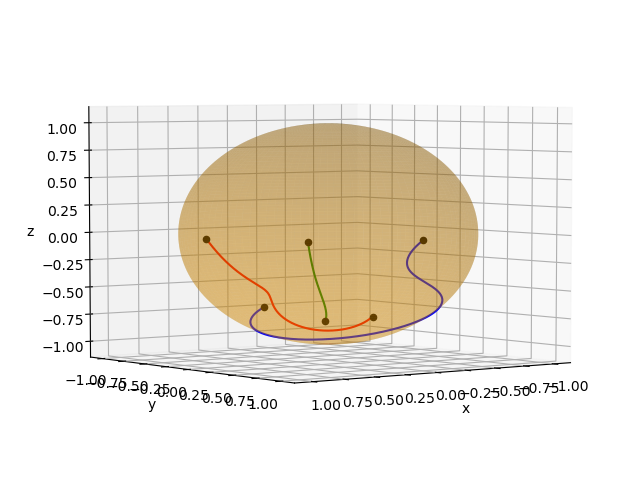}
 \includegraphics[width=6.cm]{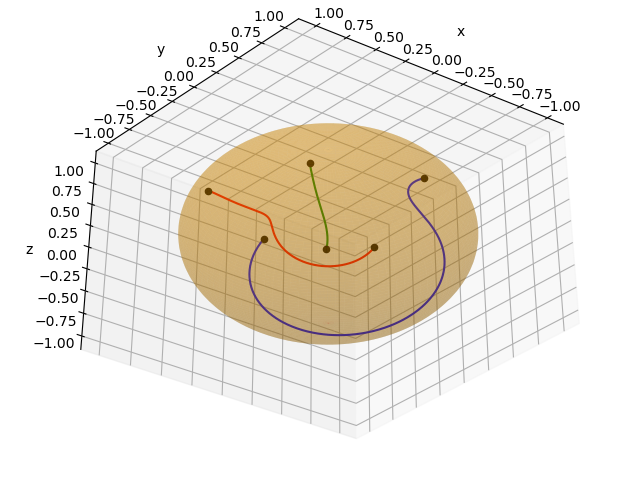}
 \caption{The figures on top show the numerical integration of the geodesics with initial conditions as defined above. We show the trajectories at time $t = 2.95$ (left)—at which point all agents fall within the collision region—and time $t = 4$ (right). The bottom figures show the solutions to the necessary conditions \eqref{eqwithexp4} satisfying the same initial conditions, from which collision avoidance is obtained.}\label{fig: s2}
 \label{figuav}
 \end{center}
\end{figure}

As in section \ref{Sec: Avoidance_bounded_domain}, we numerically calculate the constant\newline  $\displaystyle{\sum_{i=1}^s \left[(a_i^{\max})^2 + v_i^{\max}\eta_i^{\max} + v_i^0 \eta_i^0 + \frac12 \sum_{j \in \mathcal{N}_i} V_{ij}^{-} \right] < 10.32}$, from which it follows by \ref{Prop: Coll_avoid_bounded} that collision avoidance is guaranteed with tolerance $r$ if $\frac{1}{0.00001 + 256r^8} > 10.32$. This holds true for $r < 0.373$, and it can again be numerically confirmed that the solutions shown in Figure \ref{fig: s2} above avoid each other within a tolerance of $r = 0.401$.

\section*{Acknowledgements}
Both authors conduct their research at Instituto de Ciencias Matematicas
(CSIC-UAM-UC3M-UCM), Calle Nicolas Cabrera 13-15, 28049, Madrid, Spain. The project that gave rise to these results received the support of a fellowship from ”la Caixa” Foundation (ID 100010434). The fellowship codes are LCF/BQ/DI19/11730028 for Jacob R. Goodman (\textcolor{blue}{jacob.goodman@icmat.es}) and LCF/BQ/PI19/11690016 for Leonardo J. Colombo (\textcolor{blue}{leo.colombo@icmat.es}). The authors were also partially funded by Ministerio de Economia, Industria y Competitividad
(MINECO, Spain) under grant MTM2016-76702-P  and
``Severo Ochoa Programme for Centres of Excellence'' in
R$\&$D (SEV-2015-0554). All the results are original and has not been presented nor submitted to conference.


\end{document}